\documentclass[runningheads]{llncs}

\usepackage[utf8]{inputenc} 
\usepackage[T1]{fontenc}    
\usepackage{hyperref}       
\usepackage{url}            
\usepackage{booktabs}       
\usepackage{amsfonts}       
\usepackage{nicefrac}       
\usepackage{microtype}      
\usepackage{lipsum}
\usepackage[normalem]{ulem}
\usepackage{mathrsfs}
\usepackage{epsfig, graphicx}
\usepackage{latexsym,amsfonts,amsbsy,amssymb}
\usepackage{amsmath}
\usepackage{wrapfig}
\usepackage{color}
\usepackage{mathtools}
\usepackage{float}
\usepackage{cite}
\usepackage{algorithm2e}
\usepackage{amsmath}
\usepackage{amsfonts}
\usepackage{amssymb}
\usepackage{comment}
\usepackage{caption}  
\usepackage{siunitx}

\newcolumntype{L}[1]{>{\raggedright\let\newline\\\arraybackslash\hspace{0pt}}m{#1}}
\newcolumntype{C}[1]{>{\centering\let\newline\\\arraybackslash\hspace{0pt}}m{#1}}
\newcolumntype{R}[1]{>{\raggedleft\let\newline\\\arraybackslash\hspace{0pt}}m{#1}}

\RestyleAlgo{ruled}

\newcommand{\xb}{\mathbf{x}}
\newcommand{\sbb}{\mathbf{s}}
\newcommand{\gb}{\mathbf{g}}
\newcommand{\Hb}{\mathbf{H}}

\newcommand{\R}{\mathbb{R}}
\newcommand{\dd}{\delta}
\newcommand{\zero}{\mathbf{0}}
\newcommand{\trho}{\tilde\rho}
\newcommand{\ared}{\mathrm{ared}}
\newcommand{\pred}{\mathrm{pred}}
\newcommand{\ered}{\mathrm{ered}}
\newcommand{\gammainc}{\gamma_{\mathrm{inc}}}
\newcommand{\gammadec}{\gamma_{\mathrm{dec}}}
\newcommand{\true}{\mathrm{TRUE}}
\newcommand{\false}{\mathrm{FALSE}}
\usepackage{comment}

\newcommand{\algref}[1]{{\rm Algorithm~\ref{alg:#1}}}

\renewcommand{\b}[1]{\mathbf{#1}}

\usepackage{listings}
\graphicspath{{images/}}

\begin{document}

\title{TROPHY: Trust Region Optimization\\Using a Precision Hierarchy}

\author{Richard J. Clancy\inst{1,2} \and 
Matt Menickelly\inst{1} \and
Jan H\"{u}ckelheim\inst{1} \and 
Paul Hovland\inst{1} \and \\
Prani Nalluri\inst{1,3} \and
Rebecca Gjini\inst{1,4}}

\authorrunning{R.J. Clancy et al.}

\institute{Argonne National Laboratory, Lemont, IL 60439, USA \and
University of Colorado, Boulder, CO 80309, USA \and Rice University, Houston TX 77005, USA \and University of California, San Diego, CA 92093, USA}

\maketitle
\begin{abstract}
We present an algorithm to perform trust-region-based optimization for nonlinear unconstrained problems. The method selectively uses function and gradient evaluations at different floating-point precisions to reduce the overall energy consumption, storage, and communication costs; these capabilities are increasingly important in the era of exascale computing. In particular, we are motivated by a desire to improve computational efficiency for massive climate models. We employ our method on two examples: the CUTEst test set and a large-scale data assimilation problem to recover wind fields from radar returns. Although this paper is primarily a proof of concept, we show that if implemented on appropriate hardware, the use of mixed-precision can significantly reduce the computational load compared with fixed-precision solvers.

\end{abstract}


\section{Introduction}
\vspace{-1em}
\label{sec:intro}
Optimization methods are used in many applications, including engineering, science, and machine learning. The memory requirements and run time for different methods have been studied extensively and determine the problem sizes that can be run on existing hardware. Similarly, the energy consumption of each method determines its cost and carbon footprint, which is a growing concern~\cite{hao2019training}.

With the still-increasing demand for ``big data'' and ever-increasing computational power, problem scales have grown as well. To improve efficiency, modern computers tightly integrate graphical processing units (GPUs) and other accelerators. Many of these units natively support data types of differing precision to lessen the storage and computational load. Previous work has found significant differences in the overall energy consumption for double- and single-precision computations~\cite{molka2010characterizing, kestor2013quantifying}. Server-level products such as NVIDIA Tensor cores in V100 GPUs show $16 \times$ improvement over traditional double precision~\cite{abdelfattah2021survey}.

Such gains come at a cost, however. 
Classical algorithms, such as the Gram--Schmidt process, are well known to suffer from loss of orthogonality and numerical instability due to limited precision~\cite{golub1996matrix}. 
In an effort to ameliorate algorithmic issues with accuracy and stability, there has been a flurry of activity using mixed precision. These methods utilize multiple data types in a principled fashion to reduce the computational burden without sacrificing accuracy. A few of the many applications are tomographic reconstruction \cite{doucet2019mixed}, seismic modeling~\cite{ichimura2018fast}, and neural network training~\cite{jia2018highly,micikevicius2018mixed, wang2018training, carmichael2019performance}. Mixed-precision methods have been used generically to minimize the cost of linear algebra methods~\cite{abdelfattah2021survey}, iterative schemes ~\cite{strzodka2006mixed}, and improved finite element solvers ~\cite{goddeke2007performance}. 
In this paper we focus on mixed precision employed within optimization methods.

A recent paper by Gratton and Toint~\cite{gratton2020note} illustrates potential savings in an optimization setting via variable-precision trust region (TR) methods. We investigate the ideas proposed in their work but with an important difference. In particular, their algorithm (TR1DA) requires access to an approximate objective, $\bar f(\b x_k, \omega_{f,k})$, and gradient, $\bar{\b g}(\b x_k, \omega_{g,k})$, where  $\omega_{f,k} $ and $\omega_{g,k}$ are uncertainty parameters (for the $k$th iterate $\b x_k$) that satisfy
\[
    |\bar f(\b x_k, \omega_{f,k}) - f(\b x_k)| \le  \omega_{f,k}  \qquad \text{and} \qquad \frac{\| \bar{\b g}(\b x_k, \omega_{g,k}) - \b g(\b x_k) \|}{\| \bar{\b g}(\b x_k, \omega_{g,k})\|} \le \omega_{g,k}.
\]
That is, their error model requires the user to specify absolute error bounds on function and gradient values. We note that bounding error on the image becomes unlikely as computational complexity grows for reasons such as catastrophic cancellation and accumulated round-off error, making practical use of TR1DA difficult. Our focus here is on designing an algorithm that performs well without restrictive assumptions on the error resulting from lower precision.

In this paper we introduce TROPHY (\textbf{T}rust \textbf{R}egion \textbf{O}ptimization using a \textbf{P}recision \textbf{H}ierarch\textbf{Y}), a mixed-precision TR method for unconstrained optimization. We do not place restrictive assumptions on the accuracy of computed function and gradient values.
Instead, we provide practically verifiable conditions intended to determine whether the error related to a current precision level may be interfering with the dynamics of the trust region algorithm. 
If the conditions are not satisfied, we simply increase the precision level until they are. 
Our goal is to lighten the computational load without sacrificing accuracy of the final solution. By using a limited-memory, symmetric rank-1 update (L-SR1) to the approximate Hessian, the method is suitable for large-scale, high-dimensional optimization problems. We compare the method with a standard TR method---supplied with access to either a single- or double-precision evaluation of the function and gradient---on the Constrained and Unconstrained Testing Environment with safe threads (CUTEst) test problem collection~\cite{gould2015cutest} and on a large-scale weather model based on the PyDDA software package \cite{pydda2020}.

Since computational, storage, and communication savings are based on hardware implementations of different precision types rather than assumed theoretical values, our primary metric for comparison will be adjusted function evaluations rather than time. Simply put, adjusted function evaluations discount computations performed in lower-precision levels. The goal here is to provide a proof of concept for computational gains attainable by exploiting variable precision in TR methods. In practice, improvements in energy consumption, time, communication, and memory must be realized through optimized hardware, which is beyond the scope of this paper.

\section{Background} 
\label{sec:background}
\vspace{-1em}
Consider the unconstrained minimization of a differentiable function $f:\R^n\to\R$,
\begin{equation}
\label{eq:unc}
    \min_{\xb\in\R^n} f(\xb).
\end{equation}
We are motivated by problems where
the objective and its derivatives are expensive to calculate as is typical for large-scale computing. Quasi-Newton methods are a popular choice for high-dimensional problems because they empirically tend to require fewer iterations than a simple first-order method does to attain convergence.
In this paper we  focus on the TR framework, but we could have just as easily studied line-search methods such as L-BFGS, which is a popular quasi-Newton method distributed in SciPy~\cite{virtanen2020scipy}.
However, it is remarkably simpler to illustrate the effect of error on the quality of models within a TR method; 
that is likely the reason TR methods were employed in~\cite{gratton2020note}. In the following subsections we give a brief overview of the general framework for TR methods and describe the model function used in our algorithm.

\subsection{Trust Region Methods}
\vspace{-1em}
Trust region methods are iterative algorithms used for numerical optimization. At each iteration (with the counter denoted by $k$), a model function $m_k:  \R^n  \to \R$ is built around the incumbent point or iterate, $\b x_k$, such that $m_k(\b 0) = f(\b x_k)$ and $m_k(\b s) \approx f(\b x_k + \b s)$. The model, $m_k$, is intended to be a ``good'' local model of $f$ on the \emph{trust region}
$\{\sbb\in\R^n: \|\sbb\|\leq\dd_k\}$ for $\dd_k>0$.
We refer to $\dd_k$ as the \emph{trust region radius}. A \emph{trial step}, $\sbb_k$, is then computed via a(n approximate) solution to the \emph{trust region subproblem},
 \begin{equation}
  \label{eq:trsp}
  \sbb_k \  = \ \underset{\|\sbb\|\leq\dd_k}{\text{argmin}} \ \  m_k(\sbb),
\end{equation}
for $\sbb\in\R^n$.
By an approximate solution $\sbb_k$ to the TR subproblem \eqref{eq:trsp}, we mean that one requires the  \emph{Cauchy decrease condition} to be satisfied: 
\begin{equation}
    f(\b x_k) - m_k(\b s_k) \ge \frac{\mu}{2} \min \left\{ \delta_k, \frac{\| \b g_k \|}{C} \right\},
\end{equation}
where $\mu$ and $C$ are constants and $\b g_k = \b \nabla m(\b x_k)$. 
A common choice for $m_k$ is a quadratic Taylor expansion, namely, $m_k(\b s) = f(\b x_k) + \b g_k^T \b s + \frac 1 2 \b s^T \nabla^2 f(\b x_k) \b s$.
In practice, $\nabla^2 f(\b x_k)$ is typically replaced with a (quasi-Newton) approximation. 

Having computed $\sbb_k$, the standard TR method then compares the true decrease in the function value, $f(\xb_k)-f(\xb_k+\sbb_k)$, with the decrease predicted by the model,
$m_k(\zero)-m_k(\sbb_k)$. 
In particular, one computes the quantity
\begin{equation}
 \label{eq:rho_true}
 \rho_k=\displaystyle\frac{f(\xb_k)-f(\xb_k+\sbb_k)}{m_k(\zero)-m_k(\sbb_k)}.
\end{equation}
If $\rho_k$ is sufficiently positive ($\rho_k>\eta_1$ for fixed $\eta_1>0$), then the algorithm accepts $\xb_k+\sbb_k$ as the incumbent point $\xb_{k+1}$ and may possibly
increase the TR radius $\dd_{k}<\dd_{k+1}$ (if $\rho_k>\eta_2$ for fixed $\eta_2 \ge \eta_1$).
This scenario is called a \emph{successful iteration}. 
On the other hand, if $\rho_k$ is not sufficiently positive (or is negative), then the incumbent point stays the same, $\xb_{k+1}=\xb_k$, 
and we set $\dd_{k+1} < \dd_{k}$. This process is iterated until a stopping criterion is met, for instance, when the gradient norm $\|\nabla f(\xb_k)\|$ is below a given tolerance. Under mild assumptions, TR methods asymptotically converge to stationary points of $f$ \cite{conn2000trust}.

\subsection{Model Function}
\vspace{-1em}
The model function, $m_k$, must be specified for a TR algorithm. Popular choices include linear or quadratic approximations of the objective using Taylor series or interpolation methods; the latter are often employed in derivative-free optimization \cite{larson2019derivative}. Since many applications of interest are high dimensional or have costly objective and derivative functions, it is difficult, if not impossible, to compute and/or store the Hessian matrix for use in quadratic trust-region models with the memory requirement scaling as $\mathcal{O}(n^2)$. 
A common technique that exploits derivative information while keeping the cost low is to use \textit{curvature pairs} given by $\b s_k$ and $\b y_k = \nabla f(\b x_k + \b s_k) - \nabla f(\b x_k)$. After each successful iteration, the curvature pairs are used to update the current approximate Hessian denoted by $\b H_k$. These updates employ secant approximations of second derivatives. Common update rules include BFGS, DFP, and SR1 \cite{nocedal2006numerical}.

In this work we  use a limited-memory symmetric rank-1 update (L-SR1) to the approximate Hessian. This update rule requires the user to set a memory parameter that specifies a number of secant pairs to use in the approximate Hessian. Since we  require only a matrix-vector product and not the explicit Hessian, we can implement a matrix-free version reducing the storage cost to $\mathcal{O}(n)$. Thus, our TR subproblem is 
\begin{equation}
    \b s_k = \underset{\|\b s\| \le \delta_k}{\text{argmin}}  \ \ \sbb^T \nabla f(\b x_k) + \frac 1 2 \sbb ^T \b H_k \sbb,
\end{equation}
which we recast and approximately solve using the Steihaug conjugate gradient method implemented in~\cite{berahas2019quasi}[Appendix B.4].

In the next section we describe the dynamic precision framework and present criteria for when precision should switch. We then are prepared to give a formal statement of TROPHY. 
In Section~\ref{sec:testcase} we describe the problems on which we have tested TROPHY, and in Section~\ref{sec:results} we discuss the results of our experiments.

\section{Method}
\label{sec:method}
\vspace{-1em}
We suppose that we have access to a hierarchy of arithmetic precisions for the evaluation of both $f(\xb)$ and $\nabla f(\xb)$, but the direct (infinite-precision) evaluation of $f(\xb),\nabla f(\xb)$ is unavailable. We formalize this slightly by supposing we are given oracles that compute $f^p(\xb), \nabla f^p(\xb)$ for $p\in\{0,\dots,P\}$. 
With very high probability, given a uniform distribution on all possible inputs $\xb$, the $p\in\{0,\dots,P-1\}$ are such that
\begin{equation*}
 |f^p(\xb)-f(\xb)| > |f^{p+1}(\xb)-f(\xb)|, \quad \|\nabla f^p(\xb) - \nabla f(\xb)\| > \|\nabla f^{p+1}(\xb)-\nabla f(\xb)\|.
\end{equation*}
For a tangible example, if intermediate computations involved in the computation of $f(\xb)$ can be done in half, single, or double precision,
then we can denote $f^0(\xb)$ as the oracle using only half-precision computations, $f^1(\xb)$ as the oracle using only single precision, and $f^2(\xb)$ as the oracle using only double precision. 

To build on the generic TR method described in Section~\ref{sec:background}, we must specify when and how to switch precision.
We can identify two additional difficulties presented in the multiple-precision setting. 
First, it is currently unclear how to compute $\rho_k$ in \eqref{eq:rho_true}, since our error model assumes we have no access to an oracle that directly computes $f(\cdot)$.
Second, because models $m_k$ are typically constructed from function and gradient information provided by $f(\cdot)$ and $\nabla f(\cdot)$, we must identify appropriate conditions on $m_k$.

For the first of these two issues, we make a practical assumption that \textbf{the highest level of precision available to us should be treated as if it were infinite precision.}
Although this is a theoretically poor assumption, virtually all computational optimization makes this assumption implicitly; optimization methods are typically analyzed
as if computation can be done in infinite precision but are only ever implemented at some level of arithmetic precision (typically double). 
Thus, in the notation we have developed, the optimization problem we actually aim to solve is not \eqref{eq:unc} but 
\begin{equation}
    \label{eq:unc_real}
     \displaystyle\min_{\xb\in\R^n} f^P(\xb),
\end{equation}
so that the $\rho$-test in \eqref{eq:rho_true} is replaced with
\begin{equation}
     \label{eq:rho}
     \rho_k=\displaystyle\frac{f^P(\xb_k)-f^P(\xb_k+\sbb_k)}{m_k(\zero)-m_k(\sbb_k)}=\displaystyle\frac{\ared_k}{\pred_k}.
\end{equation}
We introduced $\ared$ to denote ``actual reduction'' and $\pred$ to denote ``predicted reduction.'' 
However, computing \eqref{eq:rho} still entails two evaluations of the highest-precision oracle $f^P(\cdot)$, and the intention of using mixed precision
was to avoid exactly this behavior. 
Our algorithm avoids the cost of full-precision evaluations 
by dynamically adjusting the precision level $p_k\in\{1,\dots,P\}$ between iterations so that in the $k$th iteration, $\rho_k$ is approximated by
\begin{equation}
 \label{eq:rho_observed}
 \trho_k=\displaystyle\frac{f^{p_k}(\xb_k)-f^{p_k}(\xb_k+\sbb_k)}{m_k(\zero)-m_k(\sbb_k)}=\displaystyle\frac{\ered_k}{\pred_k},
\end{equation}
introducing $\ered$ to denote ``estimated reduction.'' 
To update $p_k$, we are motivated by a strategy similar to one employed in \cite{heinkenschlossvicente} and \cite{Kouri2014InexactOF}. 
We introduce a variable $\theta_k$ that is not  initialized until the end of the first unsuccessful iteration, and we initialize $p_0=0$. When the first unsuccessful iteration is encountered, we set
\begin{equation}
    \label{eq:theta_update}
    \theta_k\gets \left|\ared_k-\ered_k\right|.
\end{equation}
Notice that we must incur the cost of two evaluations of $f^P(\cdot)$ following the first unsuccessful iteration in order to compute $\ared_k$. 
From that point on, $\theta_k$ is involved in a test triggered on every unsuccessful iteration to determine whether the precision level $p_k$ should be increased. 
Introducing a predetermined \emph{forcing sequence} $\{r_k\}$ satisfying $r_k\in[0,\infty)$ for all $k$ and $\displaystyle\lim_{k\to\infty} r_k=0$,
and fixing a parameter $\omega\in(0,1)$, 
we check on any unsuccessful iteration whether
\begin{equation}
    \label{eq:precision_test}
    \theta_k^\omega \leq \eta\min\left\{\pred_k,r_k\right\},
\end{equation}
where $\eta=\min\left\{\eta_1,1-\eta_2\right\}$. If \eqref{eq:precision_test} does not hold, then we increase $p_{k+1}=p_k+1$ and again update the value of $\theta_k$ according to \eqref{eq:theta_update} (thus incurring two more evaluations of $f^P(\cdot)$). 
The reasoning behind the test in \eqref{eq:precision_test} is that
if (the unknown) $\rho_k$ in \eqref{eq:rho} satisfies $\rho_k \geq \eta$, then
\begin{equation}\label{eq:explainer}
    \eta \leq \rho_k = \displaystyle\frac{\ared_k}{\pred_k} \leq \displaystyle\frac{|\ared_k-\ered_k|+\ered_k}{\pred_k} \approx \displaystyle\frac{\theta_k+\ered_k}{\pred_k} = \displaystyle\frac{\theta_k}{\pred_k} + \trho_k.
\end{equation}
Thus, for the practical test \eqref{eq:rho_observed} to be meaningful, we need to ensure that $\theta^k/\pred_k < \eta$, which is what \eqref{eq:precision_test} attempts to enforce. The use of $\omega$ and the forcing sequence in \eqref{eq:precision_test} is designed to ensure that we eventually do not tolerate error, since \eqref{eq:explainer} involves an approximation due to the estimate $\theta_k$. 

It remains to describe how we deal with our second identified difficulty, 
the construction of $m_k$ in the absence of evaluations of $f(\cdot)$ and $\nabla f(\cdot)$.
As is frequently done in trust  region methods, we will employ quadratic models of the form 
\begin{equation}
 \label{eq:quad_model}
 m_k(\sbb) = f_k + \gb_k^\top\sbb + \frac{1}{2}\sbb^\top\Hb_k \sbb.
\end{equation}
Having already defined rules for the update of $p_k$ through the test \eqref{eq:precision_test}, 
we take in the $k$th iteration
$f_k = f^{p_k}(\xb)$ and
$\gb_k = \nabla f^{p_k}(\xb)$.
In theory, we require $\Hb_k$ to be any Hessian approximation with a spectrum bounded above and below uniformly for all $k$.
In practice, we update $\Hb_k$ via L-SR1 updates \cite{byrd1994representations}.
By implementing a reduced-memory version, we need not store an explicit approximate Hessian, thus greatly reducing the memory cost and significantly accelerating the matrix-vector products in our model. Pseudocode for TROPHY is provided in \algref{trophy}. 

     \begin{algorithm}[h!]
    \caption{TROPHY \label{alg:trophy}}
    Initialize $0<\eta_1\leq\eta_2<1$, $\omega\in(0,1)$, $\gammainc>1$, $\gammadec\in(0,1)$, forcing seq. $\{r_k\}$.\\
    Choose initial $\dd_0>0$, $\xb_0\in\R^n$. \\
    $\theta_0\gets 0, p_k\gets 0, k\gets 0, \mathrm{failed}\gets\false$\\
    \While{some stopping criterion not satisfied}
        {
        Construct model $m_k$.\\
        (Approximately solve) \eqref{eq:trsp} to obtain $\sbb_k$.\\
        Compute $\trho_k$ as in \eqref{eq:rho_observed}.\\
        \uIf{$\trho_k>\eta_1$ \textbf{(successful iteration)}}{
        $\xb_{k+1}\gets\xb_k+\sbb_k$.\\
        \If{$\trho_k>\eta_2$ \textbf{(very successful iteration)}}{
            $\dd_{k+1}\gets \gammainc\dd_k$.\\
            }
            }
        \Else
        {
            \If{\textbf{not } $\mathrm{failed}$}{
            Compute $\theta_k$ as in \eqref{eq:theta_update}.\\
            $\mathrm{failed}\gets\true$.\\
        }
        \uIf{\eqref{eq:precision_test} holds}{
            $\dd_{k+1}\gets\gammadec\dd_k$.\\
        } 
        \Else{
            $p_{k+1}\gets p_k + 1$.\\
            Compute $\theta_k$ as in \eqref{eq:theta_update}.\\
            $\delta_{k+1}\gets\delta_k$.
        } 
        $\xb_{k+1}\gets\xb_k$. 
        }
        $k\gets k+1$.\\
    } 
\end{algorithm}
     
\section{Test Problems and Implementations}
\label{sec:testcase}
\vspace{-1em}
Our initial implementation of TROPHY is written in Python. To validate the algorithm, we focus on a well-known optimization test suite and a problem relating to climate modeling. In all cases, the algorithms terminate when one of the following conditions are met: (1) the first-order condition is satisfied, namely, $\| \nabla f^P(\b x_k) \| < \epsilon_{\text{tol}}$; (2) the TR radius is smaller than machine precision, namely, $\delta_k < \epsilon_{\text{machine}}$; or (3) the first two conditions have not been met after some maximum number of iterations. Condition 1 is a success whereas conditions 2 and 3 are failed attempts.  We describe the problem setup and implementation considerations in the current section and then discuss results in Section \ref{sec:results}.

\subsection{CUTEst}
\vspace{-1em}
For our first example we used the CUTEst set~\cite{gould2015cutest}, which is well known within the optimization community and offers a variety of problems that are challenging to solve. Each problem is given in a Standard Input Format  \cite{conn2013lancelot} file that is passed to a decoder from which Fortran subroutines are generated. The problems can be built directly by using single or double precision, making the set useful for mixed-precision comparison.

\vspace{-1em}
\subsubsection{Python Implementation of CUTEst:}
The PyCUTEst package~\cite{pycutest2018} serves as an interface between Python and CUTEst's Fortran source code. The problems are compiled via the interface; then Python scripts are generated and cached for subsequent function calls. Although CUTEst natively supports both single- and double-precision evaluations, at the time of writing, the single-precision implementations are concealed from the PyCUTEst API. 

To access single-precision evaluations, we used PREDUCER~\cite{huckelheim2019}, a Python script written to compare the effect of round-off errors in scientific computing. PREDUCER parses Fortran source code and downcasts double  data types to single. To allow for its use in existing code, all single data types are recast to double after function/gradient evaluation but before returning to the calling program. Because of the overhead associated with casting operations, we do not expect improvements in computational time. However, performance gains in terms of both accuracy and a reduction in the number of adjusted function calls for an iterative algorithm should be realized. We built the double-precision functions, too, and wrapped both functions to pass as a unified handle to TROPHY.

For the subset of unconstrained problems with dimension less than or equal to 100, we ran TROPHY with single/double switching along with the same TR method using only single or double precision. Our first-order stopping criterion was $\| \nabla f^P(\b x_k) \| < 10^{-5}$, and the maximum number of allowable iterations was 5,000. We show results for problems solved by at least one TR in Section \ref{sec:results}.

\vspace{-1em}
\subsubsection{Julia Implementation of CUTEst:}
The Julia programming language supports variable-precision floating-point data types. More precisely, it allows users to specify the number of bits used in the mantissa and expands memory for the exponent as necessary to avoid overflow. This is in contrast to the IEEE 754 standard that uses 11 (5), 24 (8), and 53 (11) bits for the mantissa (exponent) of half-, single-, and double-precision floats, respectively. In practice, one can assign enough bits for the exponent to avoid dynamic reallocation.

To exploit variable precision, we hand-coded several of the unconstrained CUTEst objectives in Julia and then computed gradients with forward-mode automatic differentiation (AD) using the ForwardDiff.jl package~\cite{RevelsLubinPapamarkou2016}. We wrote a Julia port that allows us to call  this code from Python for use in TROPHY. We opted to use forward-mode AD for its ease of implementation. In all cases used, the hand-coded Julia objective and AD gradient were compared against the Fortran implementation and found to be accurate. 

We compared TROPHY with TR methods using a fixed precision of half, single, and double precision (11, 24, and 53 bits, respectively). We used the same first-order condition of $\| \nabla f(\b x_k) \| < 10^{-5}$ but allowed this implementation to run only for 1,000 iterations. For TROPHY, we used several precision-switching sets: \{24, 53\},  \{11, 24, 53\}, \{8, 11, 17, 24, 53\}, and \{8, 13, 18, 23, 28, 33, 38, 43, 48, 53\}. The third set of precisions was motivated by the number of mantissa bits in bfloat16, fp16, fp24, fp32, and fp64, respectively. The last set increased the number of bits in increments of 5 up to double-precision.

\vspace{-1em}
\subsubsection{Multiple Doppler Radar Wind Retrieval:}
We also looked at a data assimilation problem for retrieving wind fields for convective storms from Doppler radar returns. Shapiro and Potvin \cite{shapiro2009use, potvin2012impact} proposed a method for doing so that optimizes a cost functional based on vertical vorticity, mass continuity, field smoothness, and data fidelity, among others.  Although the function calls are fairly simple, the wind field must be reconciled on a 3-D grid over space, each with an $x,\ y,$ and $z$ component. For a $39 \times 121 \times 121$ grid, there are $1,712,997$ variables. Therefore, reducing computational, storage, and communication costs where possible is paramount.  

Our work centered on the PyDDA package~\cite{pydda2020}, which was written to solve the aforementioned problem. We amended the code in two significant ways. First, to improve efficiency, we rewrote portions of the code to use JAX, an automatic differentiation package using XLA that exploits efficient computation on  GPUs \cite{jax2018github}. Since JAX natively supports single precision on  CPUs and can be recast to half and double as desired, it nicely serves as a proof of concept on a real application. Second, we modified the solver to use TROPHY rather than the SciPy implementation of L-BFGS. 

Once again, we compared TROPHY against single- and double-precision TR methods. TROPHY switched among half, single, and double precision. To avoid overflow initially for half-precision, we warm started the algorithm by providing it with the tenth iterate from the double-precision TR method, i.e., $\b x_{10}$. We perturbed this initial iterate 10 times and used the perturbed vectors as the initial guess for each algorithm (including double TR). We measured the average performance when solving each problem to different first-order conditions: $\| \nabla f(\b x_k) \| < 10^{-3}$ and $\| \nabla f(\b x_k) \| < 10^{-6}$. The maximum number of allowable iterations was 10,000.

\section{Experimental Results}
\label{sec:results}
\vspace{-1em}
We display results across the CUTEst set using data and performance profiles \cite{EDD01, JJMSMW09}. For a given metric, performance profiles help determine how a set of solvers, $\mathcal S$, performs over a set of problems, $\mathcal P$. The value $v_{ij} > 0$ denotes a particular metric (say, the final gradient norm) of the $j$th solver on problem $i$. We can then consider the performance of each solver in relation to the solver that performed best, that is, the one that achieved the smallest gradient norm. The \emph{performance ratio} is defined as
\begin{equation}
    r_{ij} = \frac{v_{ij}}{\min_j \{v_{ij}\}}.
\end{equation}
Smaller values of $r_{ij}$ are better since they are closer to optimal. The performance ratio was set to $\infty$ if the solver failed to solve the problem.
We can evaluate the performance of a solver by asking what percentage of the problems are solved within a fraction of the best. This is given by the \emph{performance profile},
\[
    h_j(\tau) = \frac{\sum_{i=1}^N \mathcal I_{\{ r_{ij} \le \tau \} }}{N},
\]
where $N = |\mathcal P|$ (the cardinality of $\mathcal P$) and $\mathcal I_{\{ A \}}$ is the indicator function such that $\mathcal I_{\{ A \}} = 1$ if $A$ is true and $0$ otherwise. Hence, better solvers have profiles that are above and to the left of the others.

\begin{figure}[t!]
    \centering
    \includegraphics[width=.9\textwidth]{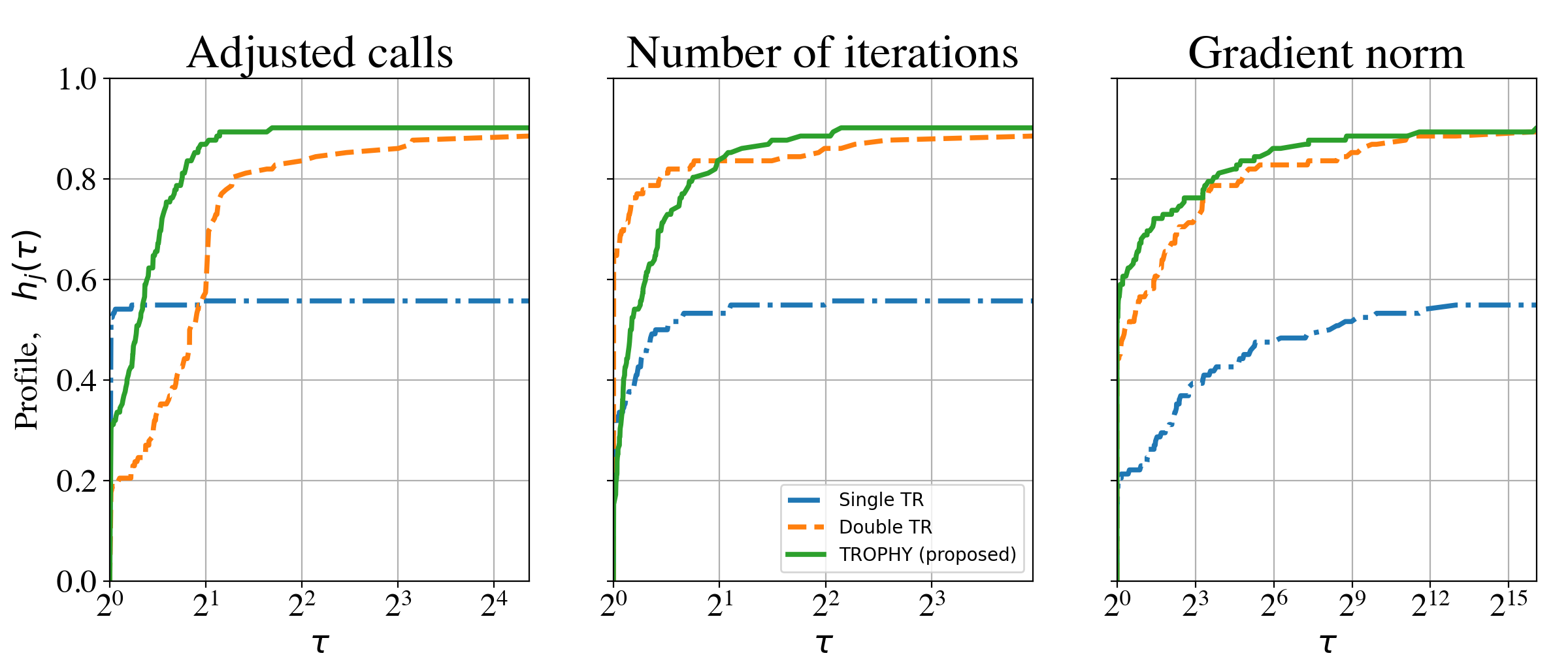}
    \caption{\small{ 
    Performance profiles for Python implementation of unconstrained CUTEst problems of dimension $< 100$ solved to first-order tolerance of $10^{-5}$. Standard single and double TR methods compared against TROPHY using single/double switching.}}
    \label{fig:pycutest_profile}
        \vspace{-.1in}
\end{figure}

\begin{figure}[t!]
    \centering
    \includegraphics[width=.9\textwidth]{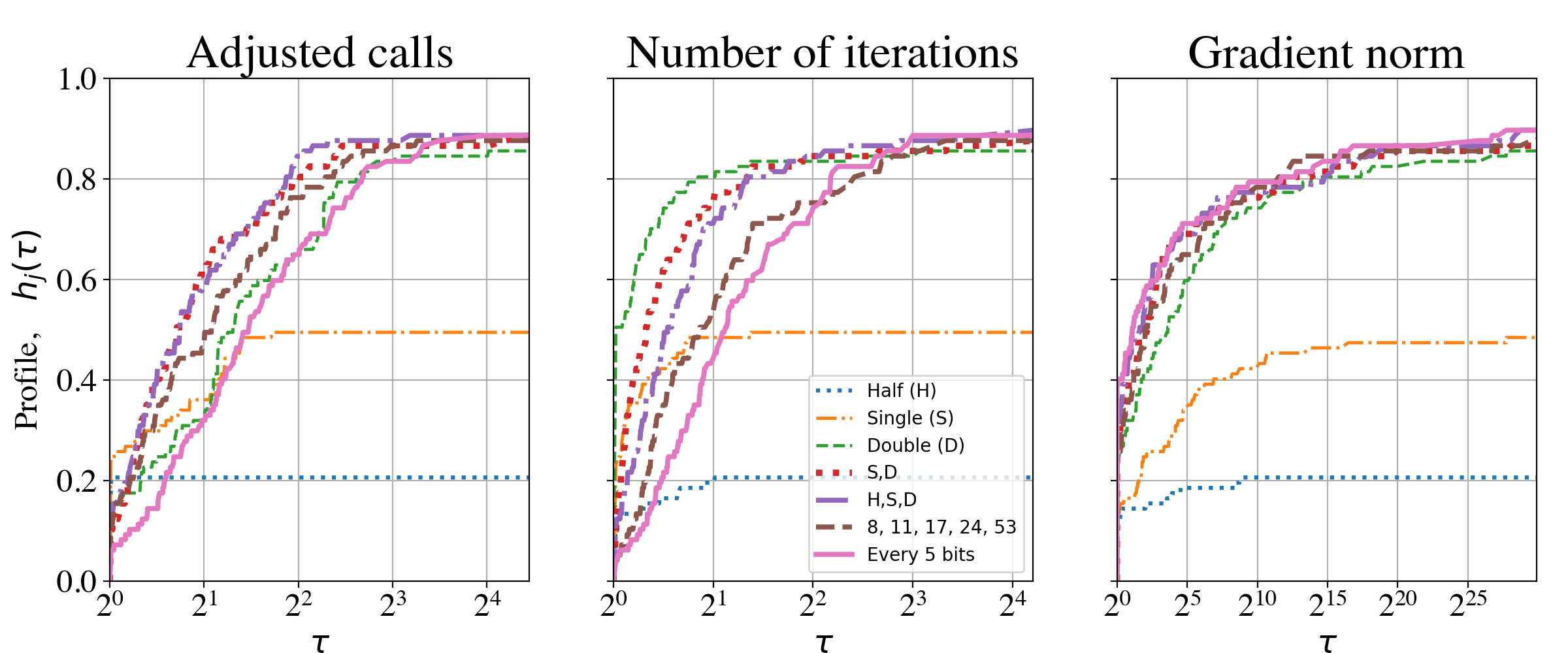}
    \caption{\small{Performance profiles for Julia implementation of unconstrained CUTEst problems of dimension $< 100$ solved to first-order tolerance of $10^{-5}$. Half, single, and double are standard TR methods using corresponding precision. ``S,D'', ``H,S,D'', ``8,11,17,24,53'', and ``Every 5 bits'' are TROPHY implementations using different precision regimes. ``Every 5 bits'' starts at 8 bits and increases to 53 bits in increments of 5 bits. A finer precision hierarchy does not imply better performance.}}
    \vspace{-.1in}
    \label{fig:julia_cutest}
\end{figure}

Motivated by the computational models in \cite{molka2010characterizing} and \cite{kestor2013quantifying}, 
we assume that the energy efficiency of single precision is between 2 and 3.6 times higher than double precision~\cite{fagan2016,galal2011}.
The storage and communication have less optimistic savings since we expect the cost of both to scale linearly with the number of bits used in the mantissa. Accordingly, we focus primarily on the model where half- and single-precision evaluations cost 1/4 and 1/2 that of a double evaluation, respectively. This gives a conservative estimate for energy cost and a favorable one for execution time. For a given problem and solver, we define \emph{adjusted calls}: 
\begin{equation} \label{eq:rel_cost}
        \text{Adj. calls} \  = \sum_{p \ \in \ \{0, 1, ..., P\}}  \frac{(\text{\# bits for prec. } p)  \times \ (\text{\# func. calls at prec. } p)}{\text{\# bits in prec. } P}.
\end{equation}
Figures \ref{fig:pycutest_profile} and \ref{fig:julia_cutest} show performance profiles for the Python and Julia implementations of CUTEst, respectively. All CUTEst problems had their first-order tolerance set to $10^{-5}$. Working from right to left in both images, we can see that the first-order condition is steady across methods provided that double-precision evaluations are ultimately available to the solver. When limited to half (11 bits) or single (24 bits), the performance suffers, and a number of problems cannot be solved. For the number of iterations in Python, we see that TROPHY and the double TR method perform comparably. The Julia implementation shows that the iterations count suffers when using low precision or TROPHY with many precision levels available for switching. This behavior is expected for low precision since the solver may never achieve the desired accuracy and hence runs longer, and for TROPHY since each precision switch requires a full iteration. For example, if 10 precision levels are available, TROPHY will take at least 10 iterations to complete. This limits the usefulness of the method on small to medium problems and problems where the initial iterate is close to the final solution.
As anticipated, TROPHY shows a distinct advantage for adjusted calls. The one exception is when there are many precision levels to cycle through, for the same reason as above. Although the initial iterate might be close to optimal, the algorithm must still visit all precision levels before breaking. The fact is made worse since each time the precision switches, two evaluations at the highest precision are required. Using two or three widely spaces precision levels yields strong results for the CUTEst set. 

\begin{small}
\begin{table}[h]
    \centering
    \caption{\small{Average performance over ten initializations for single-precision TR, double-precision TR, and TROPHY on PyDDA wind retrieval example. Adjusted calls indicate improved computational efficiency. Half, single, and double costs are 1/4, 1/2, and 1 for linear and 1/16, 1/4, and 1 for quadratic adjustments, respectively. Problem solved to $\|\nabla f(x_{\text{final}})\|_2 < 10^{-3}$ above.\\}}
     \begin{tabular}{ >{\raggedleft}p{1.9cm}| >{\centering}p{1.0cm}>{\centering}p{1.0cm} >{\centering}p{1.0cm}| >{\bfseries}>{\centering}p{1.7cm}| >{\bfseries}>{\centering}p{1.7cm}| >{\centering}p{1.5cm}  | p{1.5cm}  | }
        \hline
        Tolerance $\|\nabla f \|<10^{-3}$  &   Half calls  &  Single calls  &  Double calls &  Adj. calls (linear) &  Adj. calls (quad.) & $f_{\text{final}}$  &  $\ \|\nabla f_{\text{final}}\|$  \\
        \hline 
        \hline
        \text{Single}    &  -  & 3411  &  -  & 1706 & 853 & \num{5.3e-3}  & \num{9.5E-4}   \\
        \hline
        \text{Double}  &    -   & -  &  1877 &   1877  & 1877 & \num{4.5e-3}  &  \num{9.1e-4}    \\
        \hline
        \text{TROPHY}  &  465  &  1898 &  6  &  \uline{1071}  & \uline{510} & \num{4.7e-3}  & \num{9.4e-4} \\
        \hline
    \end{tabular}

\end{table}
\end{small}

\begin{small}
    \begin{table}[h]
    \centering
    \caption{\small{Problem solved to $\|\nabla f(x_{\text{final}})\|_2 < 10^{-6}$ accuracy. The single-precision TR method failed to converge with the TR radius falling below machine precision.}}
    \begin{tabular}{ >{\raggedleft}p{1.9cm}| >{\centering}p{1.0cm}>{\centering}p{1.0cm} >{\centering}p{1.0cm}| >{\bfseries}>{\centering}p{1.7cm}| >{\bfseries}>{\centering}p{1.7cm}| >{\centering}p{1.5cm}  | p{1.5cm}  | }
        \hline
        Tolerance $\|\nabla f \|<10^{-6}$   &   Half calls  &  Single calls  &  Double calls &  Adj. calls (linear) &  Adj. calls (quad.) & $f_{\text{final}}$  &  $\ \|\nabla f_{\text{final}}\|$  \\
        \hline 
        \hline
        \text{Single}    &  -  & $\infty$  &  -  & FAIL & FAIL & \num{9.9e-7}  & \num{3.9E-6}   \\
        \hline
        \text{Double}  &    -   & -  &  5283 &   5283  & 5283 & \num{1.8e-7}  &  \num{9.6e-7}    \\
        \hline
        \text{TROPHY}  &  465  &  7334 &  601  &  \uline{4384}  & \uline{2464} & \num{1.9e-7}  & \num{9.7e-8} \\ 
        \hline
    \end{tabular}
    \vspace{-2em}
    \label{tab:pydda}.
\end{table}
\end{small}

The wind retrieval example shows similar results. We included two ``adjusted call'' columns: one for a linear decay adjustment (memory and communication as above) and the other for quadratic decay (reduction in energy consumption). We originally iterated until $\| \nabla f(\b x_k) < 10^{-6}$ but observed that the single TR method failed to converge. Consequently, we loosened the stopping criterion 
as far as possible while maintaining the correct qualitative behavior of the solution. 
TROPHY outperformed the standard (double-precision TR) method in all cases, reducing the number of adjusted calls by 17\% to 73\%.

Our results show a promising reduction in the relative cost over naive single or double TR solvers. We expect that for many problems where function evaluations dominate linear algebra costs for the TR subproblem, our time to solve will greatly benefit from the method.

\section{Conclusion and Future Work}
\label{sec:conclusion}
\vspace{-1em}
In this paper we introduced TROPHY, a TR method that exploits variable-precision data types to lighten the computational burden of expensive function/gradient evaluations. We  illustrated proof of concept for the algorithm by implementing it on the CUTEst set and PyDDA. The full benefit of our work has not yet been realized. We look forward to implementing similar tests on hardware that can realize the full benefit of lower energy consumption and reduced memory/communication costs and ultimately shorten the time to solution. This will be especially beneficial for large scale climate models. 

We would  also like to incorporate mixed precision into line-search methods given their popularity in quasi-Newton solvers. By incorporating the same ideas into highly optimized algorithms such as the SciPy implementation of L-BFGS, we could easily deploy mixed precision to a wide population, dramatically reducing computational loads. Although TR methods are, computationally speaking, more appropriate for expensive-to-evaluate objectives, there is no reason the same ideas cannot be extended if practitioners prefer them. 

\vspace{1em}

\noindent \textbf{Acknowledgements} \\
We gratefully acknowledge the support by the Applied Mathematics activity within the U.S. Department of Energy, Office of Science, Advanced Scientific Computing Research Program, under contract number DE-AC02-06CH11357, and the computing resources provided on Swing, a high-performance computing cluster operated by the Laboratory Computing Resource Center at Argonne National Laboratory.

\small
\bibliographystyle{unsrt}  
\bibliography{references}

\scriptsize
\framebox{\parbox{4.5in}{The submitted manuscript has been created by UChicago Argonne, LLC, Operator of Argonne National Laboratory (`Argonne'). Argonne, a U.S. Department of Energy Office of Science laboratory, is operated under Contract No. DE-AC02-06CH11357. The U.S. Government retains for itself, and others acting on its behalf, a paid-up nonexclusive, irrevocable worldwide license in said article to reproduce, prepare derivative works, distribute copies to the public, and perform publicly and display publicly, by or on behalf of the Government.  The Department of Energy will provide public access to these results of federally sponsored research in accordance with the DOE Public Access Plan. \url{http://energy.gov/downloads/doe-public-access-plan}.}}
\normalsize

\end{document}